\documentclass[11pt]{amsart}
\usepackage{amssymb, amscd, amsmath}

\numberwithin{equation}{section}

\newtheorem{theorema}{Theorem}

\newtheorem{theorem}{Theorem}[section]

\newtheorem{proposition}[theorem]{Proposition}
\newtheorem{lemma}[theorem]{Lemma}

\theoremstyle{definition}
\newtheorem{definition}[theorem]{Definition}
\theoremstyle{remark}
\newtheorem{remark}[theorem]{Remark}

\newcommand\A{\mathcal{A}}

\newcommand\E{\mathcal{E}}

\newcommand{\W}{\mathcal{W}}
\newcommand{\V}{\mathcal{V}}

\renewcommand{\O}{\mathcal{O}}

\newcommand{\U}{\mathcal{U}}

\newcommand{\R}{\mathbb{R}}

\newcommand{\TT}{\mathbb{T}}

\newcommand\lie[1]{\mathfrak{#1}}

\newcommand{\fg}{\lie{g}}

\newcommand{\fk}{\lie{k}}

\newcommand{\SP}{\operatorname{Sp}}

\newcommand{\inv}{^{-1}}
\renewcommand{\ker}{ \operatorname {ker}}

\begin{document}

\title{Existence of relative periodic orbits near relative equilibria}

\author{Viktor Ginzburg and Eugene Lerman}\thanks{Supported in
part by the NSF grants
DMS-0204448 (EL) and DMS-0307484 (VG) and by the faculty
research funds of the University of California, Santa Cruz (VG)}

\address{Department of Mathematics, University of California at Santa Cruz,
Santa Cruz, CA 95064}
\address{Department of
Mathematics, University of Illinois at Urbana-Champaign, Urbana, IL 61801}
\email{ginzburg@math.ucsc.edu}
\email{lerman@math.uiuc.edu}

\date{\today}

\begin{abstract}
We show existence of relative periodic orbits (a.k.a.  relative
nonlinear normal modes) near relative equilibria of a symmetric
Hamiltonian system under an appropriate assumption on the Hessian of
the Hamiltonian.  This gives a relative version of the Moser-Weinstein
theorem.
\end{abstract}

\maketitle

\tableofcontents

\section{Introduction}

In this paper we discuss a generalization of the Weinstein-Moser 
theorem \cite{Ms,W1,W2} on the existence of  
nonlinear normal modes (i.e., periodic orbits) near an equilibrium in
a Hamiltonian system to a theorem on the existence of relative
periodic orbits (r.p.o.'s) near a relative equilibrium of a symmetric
Hamiltonian system.

More specifically let $(M, \omega _M)$ be a symplectic manifold with a proper
Hamiltonian action of a  Lie group $G$ and a corresponding equivariant
moment map $\Phi: M \to \fg^{\ast}$.    Let $h\in C^{\infty}(M)^G$ be a 
$G$-invariant
Hamiltonian.  We will refer to the quadruple $(M, \omega _M , \Phi: M \to
\fg^*, h \in C^\infty (M)^G)$ as a {\bf symmetric Hamiltonian system}.
The main result of the paper is the following theorem (the terms used
in the statement are explained below):

\begin{theorema}
	\label{main} 
Let $(M, \omega _M , \Phi: M \to \fg^*, h \in C^\infty (M)^G)$ be a
symmetric Hamiltonian system.  Suppose $x\in M$ is a positive definite
relative equilibrium of the system and let $\mu = \Phi (x)$.
Then for every sufficiently small $E>0$ the set $\{h = E + h (x) \}
\cap \Phi\inv (\mu)$ (if nonempty) contains a relative periodic
orbit of $h$.  
\end{theorema}

This theorem strengthens or complements numerous other results
about the existence of relative periodic orbits; see, e.g.,
\cite{LT,Ortega,Ortega2}.

We now define the relevant terms.  Recall that for every invariant
function $h\in C^\infty (M)^G$ the restriction $h|_{\Phi\inv (\mu)}$
descends to a continuous function $h_\mu$ on the symplectic quotient
\[
M/\!/_\mu G: = \Phi\inv (\mu)/G_\mu ,
\]
where $G_\mu $ denotes the stabilizer of $\mu \in \fg^*$ under the
coadjoint action.   If $G_\mu $ acts freely on $\Phi\inv (\mu)$ then
the symplectic quotient $M/\!/_\mu G$ is a  symplectic manifold
and the flow (of the Hamiltonian vector field) of $h$ on $\Phi\inv
(\mu)$ descends to the flow of $h_\mu$ on the quotient.  More
generally the quotient $M/\!/_\mu G$ is a symplectic stratified space:
the quotient is a union of strata, which are symplectic manifolds and
which fit together in a locally simple manner.  In this case the flow
of $h$ descends to a flow on $M/\!/_\mu G$, which is, on each stratum,
the flow of the restriction of $h_\mu$ to the stratum.\footnote{The
facts referred to above were first proved in \cite{SL} in the special
case of $G$ being compact and $\mu =0$.  The case where $G$ acts
properly and the coadjoint orbit through $\mu$ is closed was dealt
with in \cite{BL}.  The assumption on the coadjoint orbit was
subsequently removed in \cite{LW} following a suggestion in
\cite{Ortega}. }

A point $x\in \Phi \inv (\mu)$ is a {\bf relative equilibrium} of $h$
(strictly speaking of the symmetric Hamiltonian system $(M, \omega,
\Phi:M\to \fg^*, h\in C^\infty (M)^G$), if its image $[x]\in M/\!/_\mu
G$ is stationary under the flow of $h_\mu$.  If $x$ is a relative
equilibrium then there is a vector $\xi \in \fg$ such that $d h (x) =
d (\langle \Phi, \xi \rangle ) (x)$.  It is not hard to see that $\xi$
is in the Lie algebra $\fg_\mu$ of the stabilizer of $\mu$ and that
$\xi$ is unique modulo $\fg_x$, the Lie algebra of the stabilizer of
$x$. If we set
\[
h^\xi  := h - \langle \Phi, \xi \rangle ,
\]
then $d( h^\xi )(x) = 0$.  Hence the Hessian $d^2( h^\xi ) (x) $ is
well-defined.  The {\em symplectic slice} at $x\in M$ is, by
definition, the vector space
\[
V:= T_x (G\cdot x)^\omega /(T_x (G\cdot x) \cap T_x (G\cdot x)^\omega ),
\]
where $\cdot ^\omega$ denotes the symplectic perpendicular.  Note that
$V$ is naturally a symplectic representation of $G_x$, the stabilizer
of $x$.  Note also that $T_x (G\cdot x)^\omega = \ker d \Phi_x$, so
$V$ is isomorphic to a maximal symplectic subspace of $\ker d\Phi_x$.
In particular if $G_\mu$ acts freely at $x$ then $V$ models the
symplectic quotient $M/\!/_\mu G$ near $[x]$.  A computation shows
that $T_x (G\cdot x) \cap T_x (G\cdot x)^\omega = T_x(G_\mu \cdot x)$
and that $h^\xi (g\cdot x) = h^\xi (x)$ for all $g\in G_\mu$.  Hence
that subspace $T_x (G_\mu \cdot x)$ lies in the kernel of the
quadratic from $d^2 (h^\xi )(x)$.  It follows that the Hessian $d^2
(h^\xi ) (x)$ descends to a well-defined quadratic form $q$ on the
symplectic slice $V$ (which depends on $\xi$).  We say that the
relative equilibrium $x$ is {\bf positive definite} if $q$ is positive
definite for some choice of $\xi$. When the action is free, $q$ can be
thought of as the restriction of
$d^2\left(h\mid_{\Phi^{-1}(\mu)}\right)$ to the normal $V$ in
$\Phi^{-1}(\mu)$ to the orbit $G\cdot x$.  An integral curve $\gamma
(t) \subset \Phi \inv (\mu)$ of $h$ is a {\bf relative periodic orbit}
if its projection $[\gamma (t)]\subset M/\!/_\mu G$ is periodic. \\

The motivation for Theorem~\ref{main} comes from a result of Weinstein
\cite{W1} generalizing a classical theorem of Liapunov which asserts
that if $x$ is an equilibrium of a Hamiltonian $h$ on a symplectic
manifold $(M, \omega)$ and if the Hessian $d^2 h (x)$ of $h$ at $x$ is
positive definite, then for every $E>0$ sufficiently small the energy
surface
\[
\{ h = h(x) + E \} 
\]
carries at least $\frac{1}{2} \dim M$ periodic orbits. 
Now suppose $(M, \omega,
\Phi:M\to \fg^*, h\in C^\infty (M)^G$) is a symmetric Hamiltonian system, 
$\mu \in \fg^*$ is a point and the action of $G_\mu$ on $\Phi \inv
(\mu)$ is free, so that the symplectic quotient $M/\!/_\mu G$ is
smooth. If a relative equilibrium $x\in \Phi\inv (\mu)$ is positive definite
then the Hessian of $h_\mu$ at $[x]$ is positive definite.
Hence by Weinstein's theorem applied to $h_\mu$ at $[x]$,
for every $E>0$ sufficiently small the energy surface
\[
\{ h_\mu = h_\mu ([x]) + E \} 
\]
carries at least $\frac{1}{2} \dim M/\!/_\mu G$ periodic orbits.  In
other words, under these natural assumptions the manifolds
\[
\{ h = h(x) + E \} \cap \Phi\inv (\mu)
\]
carry relative periodic orbits of $h$.  It is natural to ask what
happens if $G_\mu$ does not act freely at or near the
relative equilibrium $x$ of $h$.  To address this issue let us  recall where the
strata of the symplectic quotients come from. For a subgroup $H$ of
$G$ the set $M_{(H)}$ of points of orbit type $(H)$ is defined by
\[
M_{(H)} : = \{ m\in M \mid \textrm{ the stabilizer }
 G_m \textrm{ is conjugate to } H \textrm{ in } G\}.
\]
Since the action of $G$ is proper, $M_{(H)}$ is a manifold. Moreover the set
\[
(M/\!/_\mu G)_{(H)}:= \left( M_{(H)} \cap \Phi \inv (\mu)\right)/G_\mu
\]
is naturally a symplectic manifold \cite{BL,SL}.  Now if $x\in
\Phi\inv (\mu)$ is a relative equilibrium and $G_x$ is the
stabilizer of $x$, then $(M/\!/_\mu G)_{(G_x)}$ is the stratum of the
quotient $M/\!/_\mu G$ containing $[x]$. Hence by Weinstein's theorem if
$h_\mu |_{(M/\!/_\mu G)_{(G_x)}}$ has a positive definite Hessian at $[x]$ then
the energy surface
\[
\{ h_\mu = h_\mu ([x]) + E \} 
\]
contains at least $\frac{1}{2}\dim ( M/\!/_\mu G)_{(G_x)}$ periodic orbits of 
$h_\mu |_{(M/\!/_\mu G)_{(G_x)}}$ for all $E>0$ sufficiently small.

The trivial observation above raises a natural question.  Suppose the
stratum containing the point $[x]$ is {\em not} open in the quotient
$M/\!/_\mu G$. Under suitable assumptions on the second partials of
$h$ at a relative equilibrium $x$,  must  the set
\[
\{ h_\mu = h_\mu ([x]) + E \} \subset M/\!/_\mu G
\]
contain more periodic orbits of the flow of $h_\mu$ than $\frac{1}{2}
(\dim M/\!/_\mu G)_{(H)}$?  For example, are there periodic orbits of
$h_\mu$ in nearby strata?  Theorem~\ref{main} in effect
affirmatively answer the question in a special case: the stratum of
$M/\!/_\mu G$ passing through $[x]$ is a single point $\{[x]\}$. 

Keeping in mind that $V=T_x M$, when $x$ is a fixed point of the action,
note that the following is a special case of Theorem~\ref{main} above:
\begin{theorema}\label{theorem2}
Let $K \to Sp (V, \omega)$ denote a symplectic representation of a
compact Lie group $K$ on a symplectic vector space $(V, \omega)$, let
$\Phi:V\to \fk^*$ denote the associated homogeneous moment map.  Let
$h\in C^\infty (V)^K$ be an invariant function with $dh (0) = 0$ and the
quadratic form $q:= d^2 h (0) $ positive definite.  Then for every $E>0$
sufficiently small there is a relatively periodic orbit of $h$ on the
set
\[
\{h = h(0) + E\} \cap \{ \Phi = 0\},
\]
provided the set in question is non-empty.
\end{theorema}

We will show that in fact Theorem~\ref{theorem2} implies
Theorem~\ref{main}. 

\begin{theorema}\label{theorem4}
Let $x$ be a relative equilibrium of a symmetric Hamiltonian system
$(M, \omega, \Phi: M\to \fg^*, h \in C^\infty (M)^G)$.  Denote the
stabilizer of $x$ by $G_x$, the stabilizer of $\mu = \Phi (x)$ by
$G_\mu$, the symplectic slice at $x$ by $V$ and the moment map
associated to the symplectic representation of $G_x$ on $V$ by $\Phi_V$.

Then there exists a Hamiltonian $h_V \in C^\infty (V)^{G_x}$ with 
$dh_V (0) = 0$ so that for any $E\in \R$ sufficiently small and
for any $G_x$-relatively periodic orbit of $h_V$ in
\[
\{ h_V = E\} \cap \Phi_V \inv (0)
\]
sufficiently close to $0$ there is a $G$-relatively periodic orbit of $h$ in 
\[
\{ h = h(x) +  E\} \cap \Phi\inv (\mu) .
\]
Moreover, if $x$ is a positive definite relative equilibrium of $h$,
then $h_V$ can be chosen so that the the Hessian $d^2 h_V (0) $ is
positive definite.
\end{theorema}
Clearly Theorem~\ref{theorem2} and Theorem~\ref{theorem4} together
imply Theorem~\ref{main}.
We will then reduce the proof of
Theorem~\ref{theorem2} to
\begin{theorema}\label{theorem3}
Let $Q$ be a compact manifold with a contact form $\alpha$ whose Reeb flow
generates a torus action.  Then for any contact form $\beta$
$C^2$-close to $\alpha$, the Reeb flow of $\beta$ has at least one
periodic orbit.
\end{theorema}

\subsection*{A note on notation}  

Throughout the paper the Lie algebra of a Lie group denoted by a
capital Roman letter will be denoted by the same small letter in the
fraktur font: thus $\fg$ denotes the Lie algebra of a Lie group $G$
etc.  The identity element of a Lie group is denoted by 1.  The natural 
pairing between $\fg$ and $\fg^*$ will be denoted by 
$\langle \cdot, \cdot \rangle$.

When a Lie group $G$ acts on a manifold $M$ we denote the action by an
element $g\in G$ on a point $x\in G$ by $g\cdot x$; $G\cdot x$ denotes
the $G$-orbit of $x$ and so on.  The vector field induced on $M$ by an
element $X$ of the Lie algebra $\fg$ of $G$ is denoted by $X_M$.  The
isotropy group of a point $x\in M$ is denoted by $G_x$; the Lie
algebra of $G_x$ is denoted by $\fg_x$ and is referred to as the
isotropy Lie algebra of $x$.  We recall that $\fg_x = \{ X \in \fg\mid
X_M (x) = 0\}$.  The image of a point $x\in M$ in $M/G$ under the
orbit map is denoted by $[x]$.

If $P$ is a principal $G$-bundle then $[p, m]$ denotes the point in the
associated bundle $P\times _G M = (P\times M)/G$ which is the orbit of
$(p,m) \in P\times M$.

If $\omega$ is a differential form on a manifold $M$ and $Y$ is a 
vector field on $M$, the contraction of $\omega$ by $Y$ is 
denoted by $\iota(Y) \omega$.

\section{\protect Reducing non-linear to linear: proof of 
Theorem~\ref{theorem4} }

\subsection{ Facts about  symplectic quotients}

In this subsection we gather a few facts \cite{SL,BL, LW} about
symplectic quotients that we will need in the proof of
Theorem~\ref{theorem4}.  As we mentioned in the introduction, the
symplectic quotient at $\mu \in\fg^*$ for a proper Hamiltonian action
of a Lie group $G$ on a symplectic manifold $(M, \omega)$ is the
topological space
\[
M/\!/_\mu G: = \Phi \inv (\mu) /G_\mu
\]
where as before $\Phi: M\to \fg^*$ denotes the associated equivariant
moment map.  We define the set of {\bf smooth functions} $C^\infty
(M/\!/_\mu G)$ by
\[
C^\infty (M/\!/_\mu G) = \{ f\in C^0 (M/\!/_\mu G) \mid 
\pi_\mu^* f \in C^\infty (M)^G |_{\Phi \inv (\mu)} \},
\]
where $\pi_\mu :\Phi \inv (\mu) \to M/\!/_\mu G$ denotes the orbit
map.  Note that since $\pi_\mu$ is surjective, given $h\in C^\infty
(M)^G$ there is a {\em unique} $h_\mu \in C^\infty (M/\!/_\mu G)$ with
\[
 h|_{\Phi\inv (\mu) } = \pi_\mu ^* h_\mu .
\]
One often refers to $h_\mu$ as the {\bf reduction} of the Hamiltonian $h$ at $\mu$.

\begin{theorem}[Arms-Cushman-Gotay \cite{acg}] \label{fact1}
The Poisson bracket $\{ \cdot , \cdot \}$ on $C^\infty (M)$ induces a
Poisson bracket $\{ \cdot , \cdot \}_\mu$ on $C^\infty (M/\!/_\mu G)$
so that
\[
\pi_\mu ^* : C^\infty (M/\!/_\mu G) \to C^\infty (M)^G |_{\Phi \inv (\mu)}
\]
is a Poisson map.
\end{theorem}

\begin{definition}
We say that a curve $\gamma : I \to M/\!/_\mu G$ is {\bf smooth}
($C^\infty$) if $f\circ \gamma : I \to \R$ is smooth for any $f \in
C^\infty (M/\!/_\mu G)$ ($I$, of course, is an interval).  A curve
$\gamma : I \to M/\!/_\mu G$ is an {\bf integral curve} of a function
$f\in C^\infty (M/\!/_\mu G)$ if for any $k\in C^\infty (M/\!/_\mu G)$
we have
\[
\frac{d}{dt} (k\circ \gamma ) (t)  = \{ f, k\} _\mu (\gamma (t)).
\]
\end{definition}
The following fact is an easy consequence of the well-known result 
that for proper group
actions (the only kind of actions we consider) smooth invariant functions
separate orbits:

\begin{proposition}\label{fact2}
Integral curves of functions in $C^\infty (M/\!/_\mu G)$ are unique.
\end{proposition}

\noindent
It is also not hard to see that Theorem~\ref{fact1} implies that if
$\gamma$ is an integral curve of an invariant Hamiltonian $h\in
C^\infty (M)^G$ lying in $\Phi\inv (\mu)$ then $\pi_\mu \circ \gamma$
is an integral curve of the corresponding reduced Hamiltonian $h_\mu
\in C ^\infty (M/\!/_\mu G)$ (as above $\pi_\mu ^* h_\mu = h |_{\Phi
\inv (\mu)}$).  Combining this with Proposition~\ref{fact2} we get
\begin{lemma}
Let $h\in C^\infty (M)^G$ be an invariant Hamiltonian and $h_\mu\in
C^\infty (M/\!/_\mu G)$ the corresponding reduced Hamiltonian at
$\mu$.  If $\gamma :I \to M/\!/_\mu G$ is an integral curve of $h_\mu$
there exists an integral curve $\tilde{\gamma} : I \to \Phi\inv (\mu)$
of $h$ so that
\[
\pi _\mu \circ \tilde{\gamma}  = \gamma .
\]
\end{lemma} 
Note that $\tilde{\gamma}$ is not unique: for any $a\in G_\mu$ the
curve $a\cdot \tilde{\gamma}$ is also an integral curve of $\gamma$
projecting down to $\gamma$.  We will need one more fact about
integral curves of functions in symplectic quotients, which is an easy
consequence of Proposition~\ref{fact2}.

\begin{lemma}\label{lemma2.3}
Suppose $\tau: M/\!/_\mu G \to M' /\!/_{\mu'} G'$ is a continuous map
between two symplectic quotients such that the pull-back $\tau^* $
maps $C^\infty (M' /\!/_{\mu'} G')$ to $C^\infty ( M/\!/_\mu G)$
preserving the Poisson brackets (i.e, $\tau $ is a morphism of
symplectic quotients).  Then for any $h' \in C^\infty (M' /\!/_{\mu'}
G')$ if $\gamma$ is an integral curve of $\tau^* h'$ then $\tau \circ
\gamma$ is an integral curve $h'$.
\end{lemma}

\begin{proof}
Since $\gamma$ is an integral curve of $\tau^* h'$ 
\[
\frac{d}{dt} ((\tau^* f') \circ \gamma )(t) = 
\{ \tau^* h', \tau^* f' \} _\mu (\gamma (t)).
\]
for any $f' \in C^\infty (M' /\!/_{\mu'} G')$.  Hence
\[
\frac{d}{dt} f' (\tau \circ \gamma) =  \{ \tau^* h', \tau^* f' \} _\mu (\gamma )
=   \tau^* (\{ h' , f' \}_{\mu'} )(\gamma)=
\{ h' , f'\}_{\mu'} (\tau \circ \gamma).
\]
\end{proof}
\noindent
This ends our digression on the subject of symplectic quotients.\\

In proving Theorem~\ref{theorem4} we will argue that there is a
$G_x$-equivariant symplectic embedding $\sigma:\V \to \U$ of a
$G_x$-invariant neighborhood $\V$ of $0$ in $V$ into a $G$-invariant
neighborhood $\U$ of $x$ in $M$ which induces a morphism
\[
\bar{\sigma} : \V/\!/G_x \to \U/\!/_\mu G
\]
of symplectic quotients so that $\bar{\sigma}$ embeds
$\V/\!/G_x$ as a connected component of $\U/\!/_\mu G$.  Note that for
$\sigma$ to induce $\bar{\sigma}$ we would want $\sigma$ to map
$\Phi_V \inv (0) \cap \V$ into $\Phi \inv (\mu) \cap \U$ in such a way
that the diagram
\begin{equation}
\label{star}
\begin{CD}
{\Phi_V\inv(0) \cap \V} @> {\sigma} >> {\Phi\inv (\mu) \cap \U} \\
@ V{\pi_0} VV @ VV {\pi_\mu} V \\
{\V/\!/_0 G_x} @> {\bar{\sigma}}> > {\U/\!/_\mu G} \\
\end{CD}
\end{equation}
commutes, where $\pi_0$, $\pi_\mu$ are the respective orbit maps.  As
before given $f\in C^\infty (\V)^{G_x}$ we denote by $f_0$ the
unique function in $C^\infty (\V/\!/_0 G_x)$ with $f|_{\Phi _V \inv
(0)\cap \V} = \pi_0^* f_0$ and similarly $h_\mu \in C^\infty
(\U/\!/_\mu G)$ is determined by $\pi^*_\mu h_\mu= h|_{\Phi^{-1} (\mu)}$.  
Then the commutativity of (\ref{star}) implies:
\[
(\sigma^* h)_0 = \bar{\sigma}^* h_\mu
\]
for any $h\in C^\infty (\U)^G$.

Suppose next that we have constructed $\sigma: \V \to \U$ with the desired properties.  Given $h\in C^\infty (M)^G$ and any $\xi \in \fg_\mu$ let 
\[
h_V = \sigma ^* (h - \langle \Phi , \xi \rangle ).
\]
Then $h_V  \in C^\infty (\V)^{G_x}$.  Moreover, since (\ref{star}) commutes,
\[
(h_V)_0 = (\sigma ^* (h - \langle \Phi , \xi \rangle ))_0 = 
\bar{\sigma}^* (h - \langle \Phi , \xi \rangle )_\mu  = \bar{\sigma }^* h_\mu - \langle \mu, \xi \rangle.
\]
In other words,
\begin{equation}\label{**}
(h_V)_0 = \bar{\sigma}^* h_\mu  + \textrm{ constant}.
\end{equation}
Lemma~\ref{lemma2.3} and equation~(\ref{**}) imply that if $\gamma_V$
is an integral curve of $h_V$ in $\V \cap \Phi _V \inv (0) \cap \{h_V
= h_V (0) + E\} $ then $\pi _0 \circ \gamma _V$ is an integral curve
of $(h_V)_0$ in $\{(h_V)_0 = (h_V) (\pi_0(0)) + E\} $.  Hence
$\bar{\sigma} \circ \pi_0 \circ \gamma _V$ is an integral curve of
$h_\mu$ in $\{ h_\mu = h_\mu (\pi _\mu (x)) + E\}$.  It follows that
there is an integral curve $\gamma$ of $h$ in $\{ h = h (x) + E\} \cap
\U \cap \Phi \inv (\mu)$ with
\[
\pi _\mu \circ \gamma = \bar{\sigma } \circ \pi_0 \circ \gamma_V.
\]
If $\gamma _V$ is a $G_x$-relative periodic orbit of $h_V$ then $\pi_0\circ \gamma_V $ is a periodic orbit of $(h_V)_0$.  Consequently $\gamma$ is a $G$-relative periodic orbit 
of $h$.

Finally note that if additionally we can arrange for 
\[
d\sigma _0 (T_0 \V) \subset T_x (G\cdot x)^\omega,
\]
then since $\sigma$ is symplectic 
\[
d\sigma _0 (T_0 \V) \cap T_x (G_\mu \cdot x)  = \{0\}. 
\]
Consequently if $\xi \in \fg_\mu$ is such that $d (h - \langle \Phi, \xi \rangle) (x) = 0$ and 
$d^2 (h - \langle \Phi, \xi \rangle) (x)|_{T_x (G\cdot x)^\omega }$ is positive semi-definite of maximal rank, then
\[
d^2 (h - \langle \Phi, \xi \rangle) (x)|_{d\sigma_0 (T_0 \V)}
\]
is positive definite.  Therefore the Hessian
\[
d^2 (\sigma ^* (h - \langle \Phi, \xi \rangle) ) (0)
\]
is positive definite as well.  We conclude that in order to prove Theorem~\ref{theorem4}
it is enough to construct the embedding $\sigma$ with the desired properties.  In other words it is enough to prove:

\begin{proposition} \label{prop1}
Let $(M, \omega, \Phi: M\to \fg^*) $ be a symplectic manifold with a proper Hamiltonian action of a Lie group $G$.  Fix a point $x$ in $M$.  Let $\mu = \Phi (x)$, let $\Phi_V: V\to \fg_x^*$ be the homogeneous moment map associated with the symplectic slice representation $G_x  \to \SP (V, \omega_V)$.  There exists a $G_x$-invariant neighborhood $\V$ of $0$ in $V$, a $G$-invariant neighborhood $\U$ of $x$ in $M$ and a
$G_x$-equivariant embedding $\sigma :\V \to \U$ with
\[
\sigma (\Phi_V \inv (0) \cap \V) \subset \Phi\inv (\mu) \cap \U
\]
such that the composition
\[
\Phi _V \inv (0) \cap \V \stackrel{\sigma}{\to} \Phi \inv (\mu) \cap \U \stackrel{\pi_\mu}{\to}
(\Phi\inv (\mu) \cap \U)/G_\mu = \U/\!/_\mu G
\]
drops  down to 
\[
\bar{\sigma} : \V/\!/_0 G_x = (\Phi_V \inv (0)  \cap \V) / G_x \to \U/\!/ _\mu G
\]
making (\ref{star}) commute.  Moreover,

\begin{enumerate}
  \item $\bar{\sigma}(\V/\!/_0 G_x )$ is a connected component of $\U/\!/ _\mu G$ and $\bar{\sigma}$ is a homeomorphism onto its image;
  \item the pull-back $\bar{\sigma}^*$ sends  $C^\infty (\U/\!/ _\mu G)$ isomorphically to 
  $C^\infty (\V/\!/_0 G_x)$ (as Poisson algebras).
\end{enumerate}
Additionally
\[
d\sigma_0 (T_0 \V) \subset T_x (G\cdot x)^\omega .
\]
\end{proposition}
\noindent
Proposition~\ref{prop1} will follow from 

\begin{proposition} \label{prop2}
As above let $(M, \omega, \Phi: M\to \fg^*) $ be a symplectic manifold
with a proper Hamiltonian action of a Lie group $G$.  Fix a point $x$
in $M$.  Let $\mu = \Phi (x)$, let $\Phi_V: V\to \fg_x^*$ be the
homogeneous moment map associated with the symplectic slice
representation $G_x \to \SP (V, \omega_V)$.  There exists a slice
$\Sigma$ at $x$ for the action of $G$ on $M$, a $G_x$-invariant
neighborhood $\V$ of $0$ in $V$ and a $G_x$-equivariant embedding
$\sigma :\V \to \Sigma$ so that
\begin{enumerate}
  \item $\sigma (\V) $ is closed in $\Sigma$;
  \item $\sigma (\Phi_V\inv (0) \cap \V) = \Sigma \cap \Phi\inv (\mu)$;
  \item $\sigma^* \omega = \omega _V$ and $\Phi \circ \sigma = i\circ \Phi_V + \mu$
  where $i: \fg_x^* \to \fg^*$ is a $G_x$-equivariant injection;
\item $G_\mu \cdot (\Sigma \cap \Phi\inv (\mu)) $ is a connected component of 
$\Phi \inv (\mu) \cap \U$, where $\U = G\cdot \Sigma$.
\end{enumerate}
\end{proposition}

\begin{proof}[Proof of  Proposition~\ref{prop2}]
Our proof of Proposition~\ref{prop2} uses the Bates-Lerman version
\cite{BL}[pp. 212--215] of the local normal form theorem for moment
maps of Marle, Guillemin and Sternberg:

\begin{theorem}[\cite{BL}] \label{symplectic_model}
Let $(M, \omega)$ be a symplectic manifold with a proper Hamiltonian
action of a Lie group $G$ and a corresponding equivariant moment map
$\Phi: M\to \fg^* $. Fix $x\in M$, let $\mu = \Phi (x)$, $(V,
\omega_V)$ the symplectic  slice at $x$, $\Phi_V :V \to \fg_x^*$ the associated 
homogeneous moment map.  Choose a $G_x$-equivariant splitting
\[
\fg^* = \fg_x ^* \oplus (\fg_\mu/\fg_x)^* \oplus \fg_\mu ^\circ
\]
($\fg_\mu^\circ$ denotes the annihilator of $\fg_\mu$ in $\fg^*$) and
thereby $G_x$-equivariant injections
\[
i: \fg_x^* \to \fg^*, \quad j:(\fg_\mu/\fg_x)^* \to \fg^*.
\]
Let 
\[ 
Y = G\times _{G_x} ((\fg_\mu/\fg_x)^* \times V);
\]
it is a homogeneous vector bundle over $G/G_x$.  There exists a closed
2-form $\omega_Y$ on $Y$ which is non-degenerate in a neighborhood of the
zero section $G\times _{G_x} (\{(0, 0)|\}) = G\cdot [1, 0, 0]$ (1
denotes the identity in $G$) such that
\begin{enumerate}
\item  a $G$-invariant neighborhood of $G\cdot x$ in $(M, \omega)$ is 
$G$-equivariantly symplectomorphic to a neighborhood of $G\cdot [1, 0, 0]$ 
in $(Y, \omega_Y)$;
\item the moment map  $\Phi_Y$ for the action of $G$ on $(Y, \omega_Y)$ is given by
\[
 \Phi_Y ([g, \eta, v])  = g\cdot (\mu + j(\eta) + i (\Phi_V (x)))
\]
for all $(g, \eta, v) \in G\times (\fg_\mu/\fg_x)^* \times V$;
\item the embedding $\iota :V \to Y$, $\iota (v) = [1, 0, v]$ is symplectic:
$\iota^* \omega_Y  = \omega_V$.
\end{enumerate}
\end{theorem}
\noindent
Hence, we can assume without loss of generality that $(M,
\omega, \Phi) = (Y, \omega_Y, \Phi_Y)$ and $x = [1, 0, 0]$.  Note that the 
embedding $\kappa: (\fg_\mu/\fg_x)^* \times V \to Y$, $\kappa (\eta, v) = [1,
\eta, v]$ is a slice at $x$ for the action of $G$ on $Y$.

We now argue that for a small enough $G_x$-invariant neighborhood $\V$
of $0$ in $V$ and a small enough $G_x$-invariant neighborhood $\W$ of
$0$ in $(\fg_\mu/\fg_x)^*$
\[
\Sigma := \kappa (\W \times \V) \subset Y
\]
is the desired slice and 
\[
\sigma = \kappa |_{\{0\} \times \V} : \V \to \Sigma, \quad \sigma (v) = [1, 0, v]
\]
is the desired embedding.  

Note that no matter how $\V$ and $\W$ are chosen we automatically
have that $\sigma(\V)$ is closed in $\Sigma$ and $\sigma^* \omega_Y =
\omega_V$.  Hence $\sigma(\V)$ is closed in $\U := G\cdot \Sigma$,
which is a $G$-invariant neighborhood of $x$. Note also that 
\[
\Phi_Y \circ \sigma (v) = \Phi _Y ([1, 0, v]) = \mu + i (\Phi _V (v)).
\]
Next we make our choice of $\V$ and $\W$ and prove that the resulting
embedding $\sigma$ has all the desired properties. For this purpose,
factor $\Phi_Y : Y \to \fg^*$ as a sequence of maps (we identify
$\fg_\mu^*$ with $j((\fg_\mu/\fg_x)^*) \oplus i(\fg_x^*) \subset
\fg^*$):
\[
G\times _{G_x} ((\fg_\mu/\fg_x)^* \times V) \stackrel{F_1}{\to} 
G\times _{G_x} ((\fg_\mu/\fg_x)^* \times \fg_x^*)\stackrel{F_2}{\to} 
G\times _{G_\mu} \fg_\mu^* \stackrel{\E}{\to } \fg^*,
\]
where
\begin{align*}
F_1 ([g, \eta, v]) &= [g, \eta, \Phi_V (v)],\\
F_2 ([g, \eta, \nu])&= [g, j(\eta) + i(\nu)] \textrm{ and}\\
\E ([g, \vartheta] )&= g\cdot (\mu + \vartheta) .
\end{align*}
Since the tangent space $T_\mu (G\cdot \mu)$ is canonically isomorphic
to the annihilator $\fg_\mu^\circ$ and since $\fg^* = \fg_\mu ^*
\oplus \fg_\mu$, the vector bundle 
$G \times _{G_\mu} \fg_\mu^*$ is the normal bundle for the embedding
$G\cdot \mu \hookrightarrow \fg^*$ and $\E: G \times _{G_\mu}
\fg_\mu^* \to \fg^*$ is the exponential map for a flat $G_x$-invariant
metric on $\fg_x^*$.  Therefore $\E$ is a local diffeomorphism near
the zero section.  In particular there is a $G_x$-invariant
neighborhood $\O$ of $[1,0]\in G\times _{G_\mu} \fg_\mu^*$ so that
$\E|_\O$ is a diffeomorphism onto its image.  Let $\O' = (F_2 \circ
F_1) \inv (\O)$.  Then
\begin{align*}
\O' \cap \Phi_Y \inv (\mu)  &= \O' \cap (F_2 \circ
F_1) \inv ((\E|_\O)\inv (\mu)) \\
&= \O' \cap F_1\inv (F_2\inv ([1, 0]))\\
&= \O' \cap F_1 \inv (G_\mu \times _{G_x} \{ (0,0)\})\\
&= \O' \cap G_\mu \times _{G_x} (\{0\}  \times \Phi_V \inv (0)).
\end{align*}
We may take $\O$ to be of the form $\A \times_{G_\mu} (\W \times \V')$
where $\A\subset G$ is a $G_x \times G_\mu$-invariant neighborhood of
1, $\W \subset (\fg_\mu/\fg_x)^*$ is a $G_x$-invariant neighborhood of
$0$ and $\V' \subset \fg_x^*$ is a convex $G_x$-invariant neighborhood
of $0$.  We take $\V = \Phi_V \inv (\V').$ Then $\V \cap \Phi_V
\inv (0)$ is connected.  With the choices above,
$\O' = \A \times _{G_x} (\W \times \V)$ and 
\[ 
\O'  \cap \Phi_Y \inv (\mu) = \left(\A \times
_{G_x} (\W \times \V)\right)  \cap \left( G_\mu \times _{G_x} (\{0\} \times (\Phi_V
\inv (0) \cap \V))\right).
\]
Since $ G_\mu \times _{G_x} (\{0\} \times (\Phi_V \inv (0) \cap \V))$
is closed in $G\times _{G_x} (\W \times
\V) = G\cdot \Sigma = \U$ and since $\Phi_V \inv (0) \cap \V$ is
connected, $G_\mu \cdot \sigma (\Phi _V \inv (0) \cap \V) = G_\mu
\times _{G_x} (\{0\} \times (\Phi_V \inv (0) \cap
\V))$ is a connected component of $\Phi_Y\inv (\mu) \cap \U$.  This proves 
property~(4).

Note that $\Sigma = \{ [1, \eta, v] \mid \eta \in \W, v\in \V\} \subset \A
\times _{G_x} (\W \times \V)$.  Hence $\Phi_Y \inv (\mu) \cap \Sigma =
(\Phi_Y \inv (\mu) \cap \O' ) \cap \Sigma = \{ [1, 0, v] \mid v\in
\Phi_V\inv (0) \cap \V\}$, i.e.,
\[
\sigma (\Phi_V \inv (0) \cap \V) = \Phi_Y \inv (\mu) \cap \Sigma,
\]
which proves property~(2) and thereby finishes the proof of
Proposition~\ref{prop2}.
\end{proof}

\begin{proof}[Proof of Proposition~\ref{prop1}]
We continue to use the notation above.  Since $\sigma (\Phi_V \inv (0)
\cap \V) = \Phi \inv (\mu) \cap \Sigma$ and since $\sigma: \V \to
\Sigma$ is a closed embedding, the restriction $\sigma |_{\Phi_V \inv
(0) \cap \V} :\Phi_V \inv (0) \cap \V\to \Phi\inv (\mu ) \cap \Sigma$
is a $G_x$-equivariant homeomorphism.  Hence 
\[ 
\bar{\sigma} : (\Phi_V\inv (0) \cap \V)/G_x \to 
(\Phi\inv (\mu ) \cap \Sigma )/G_x 
\] 
is a homeomorphism as well.  Since $\Sigma$ is a slice at $x$ for the
action of $G$ on $M$, it is also a slice for the action of $G_\mu$.
Consequently
\[
(G_\mu \cdot (\Phi\inv (\mu ) \cap \Sigma ))/G_\mu \cong
 (\Phi\inv (\mu ) \cap \Sigma )/G_x .
\]
Since $G_\mu \cdot (\Phi\inv (\mu ) \cap \Sigma )$ is a component of
$\Phi\inv (\mu) \cap \U$, $\bar{\sigma} : \V/\!/_0 G_x \to \U/\!/_\mu
G$ is a homeomorphism onto its image.  Moreover, the diagram
(\ref{star}) commutes.

We now argue that $\bar{\sigma}$ pulls back the smooth functions in
$C^\infty (\V/\!/G_x)$ to smooth functions in $C^\infty (\U/\!/_\mu
G)$ and that the pull-back is an isomorphism of Poisson algebras.
Since $\Sigma$ is a slice and $\U = G \cdot \Sigma$, the restriction
\[
C^\infty (\U)^G \to C^\infty (\Sigma)^{G_x}, \quad f\mapsto f|_\Sigma
\]
is a bijection.    Since $\sigma (\Phi_V\inv(0) \cap \V) \subset 
\Phi\inv (\mu) \cap \U$ and since 
\[
\begin{CD}
{C^\infty (\U)^G} @> {\sigma^*} >> C^\infty (\V)^{G_x}\\ @ V{} VV @ VV {} V \\ 
{C^\infty(\U)^G |_{\Phi\inv (\mu) \cap \U} }@> {\bar{\sigma}^*}> > 
{C^\infty (\V)^{G_x}|_{\Phi_V \inv (0) \cap \V}} \\
\end{CD}
\]
commutes, where the vertical arrows are restrictions, and since the
top arrow is surjective, the bottom arrow is surjective as well.
Since $\Sigma$ is a slice and $\U = G\cdot \Sigma$, any function $f\in
C^\infty (\U)^G |_{\Phi\inv (\mu) \cap \U}$ is uniquely defined by its
values on $\Phi\inv (\mu) \cap \Sigma = \sigma (\Phi_V \inv (0) \cap
\V)$.  Hence the bottom arrow $\bar{\sigma}^*$ is also
injective. Since $\sigma: \V \to \U$ is symplectic, $\sigma^* :
C^\infty (\U) \to C^\infty (\V)$ is Poisson.  Hence $\sigma ^*:
C^\infty (\U)^G \to C^\infty (\V)^{G_x}$ is also Poisson.
Consequently $\sigma^*: C^\infty(\U)^G |_{\Phi\inv (\mu) \cap \U }\to
C^\infty (\V)^{G_x} |_{\Phi_V \inv (0) \cap \V}$ is Poisson as well.
We conclude that
\[
\bar{\sigma}^* : C^\infty (\U/\!/_\mu G) \to C^\infty (\V /\!/_0 G_x)
\]
is an isomorphism of Poisson algebras.

Finally since $\Phi \circ \sigma = i\circ \Phi _V + \mu$,
\[
d \Phi _x \circ d\sigma _0 = i \circ d( \Phi_V) _0 .
\]
Since $\Phi_V$ is quadratic homogeneous, $ d( \Phi_V) _0 =0 $.  Hence 
\[
d\sigma _0 (T_0 \V) \subset \ker d\Phi _x = T_x (G\cdot x)^\omega .
\]
\end{proof}
This concludes our proof of Theorem~\ref{theorem4} as well.

\section{From invariant Hamiltonians on vector spaces to Reeb flows: 
Theorem~\ref{theorem3} implies Theorem~\ref{theorem2}}

\begin{remark}
Theorem~\ref{theorem2} is easily seen to be true in a special case:
the set $V^K$ of $K$-fixed vectors is a subspace of $V$ of positive
dimension.  Indeed, since $h$ is $K$-invariant its Hamiltonian flow
preserves the symplectic subspace $V^K$, which is contained in the
zero level set $\Phi\inv (0)$ of the moment map.  Moreover the flow of
$h$ in $V^K$ is the Hamiltonian flow of the restriction $h|_{V^K}$.
Hence Weinstein's theorem applied to the Hamiltonian system $(V^K,
\omega |_{V^K}, h|_{V^K})$ guarantees that for any $E>0$ sufficiently
small there are at least $\frac{1}{2} \dim V^K$ periodic orbits of
$h|_{V^K}$ in the surface
\[
\{h|_{V^K} =E\} = \{h =E \} \cap V^K \subset  \{h =E \} \cap \{ \Phi = 0\}.
\]
\end{remark}
To show that Theorem~\ref{theorem3} implies Theorem~\ref{theorem2} we
first need to digress on the subject of contact quotients.

\subsection{Facts about contact quotients}

Suppose that a Lie group $G$ acts properly on a manifold $\Sigma$,
preserving a contact form $\beta$.  The associated moment map $\Psi
:\Sigma \to \fg^*$ is defined by
\[
\langle \Psi (x) , \xi \rangle = \beta _x (\xi _M (x))
\]
for all $x\in \Sigma$, all $\xi \in \fg^*$.  The map $\Psi $ is
$G$-equivariant.  The contact quotient at zero $\Sigma /\!/ G$ is, by
definition, the set
\[
\Sigma /\!/ G := \Psi \inv (0)/G .
\]
Just as in the case of symplectic quotients the contact quotients are
stratified spaces \cite{LW} with the stratification induced by the
orbit type decomposition:
\[
\Sigma /\!/ G := \coprod  _{H< G} (\Psi \inv (0) \cap \Sigma _{(H)}) /G ,
\]
where the disjoint union is taken over conjugacy classes of subgroups of $G$.
Additionally each stratum
\[
(\Sigma /\!/ G )_{(H)}:=  (\Psi \inv (0) \cap \Sigma _{(H)}) /G
\]
is a contact manifold and the contact form $\beta _{(H)}$ on each
stratum is induced by the contact form $\beta $ on $\Sigma$
\cite[Theorem~3, p. 4256]{Wi}.  More precisely for each subgroup $H$
of $G$ the set $\Psi \inv (0) \cap \Sigma _{(H)}$ is a manifold and
\[
\pi_{(H)} ^* \beta _{(H)} = \beta |_{\Psi \inv (0) \cap \Sigma _{(H)}}, 
\]
where $\pi_{(H)} : \Psi \inv (0) \cap \Sigma _{(H)} \to (\Psi \inv (0)
\cap \Sigma _{(H)}) / G = (\Sigma /\!/ G )_{(H)}$ is the orbit map.
It is not hard to see that the flow of the Reeb vector field $X$ of
$\beta$ preserves the moment map and the orbit type decomposition,
hence descends to a strata-preserving flow on the quotient $\Sigma
/\!/ G$.  Also, on each stratum the induced flow is the Reeb flow of
the induced contact form $\beta _{(H)} $.\\

We're now ready to prove that Theorem~\ref{theorem3} implies
Theorem~\ref{theorem2}.  It is no loss of generality to assume that
$h(0) = 0$. Since the quadratic form $q = d^2 h (0)$ is positive
definite, the energy surface
\[
\{q = E\},
\]
$E>0$, is a $K$-invariant hypersurface star-shaped about 0.  Hence 
\[
\alpha _E := \iota (R) \omega |_{\{q = E\}}
\]
is a $K$-invariant contact form, where $R (v) = v$ denotes the radial
vector field on $V$.  For $E>0$ sufficiently small, the $K$-invariant
set
\[
\{ h = E\}
\] 
is a hypersurface which is $C^2$-close to $\{q = E\}$. Hence 
\[
\beta  _E := \iota (R) \omega |_{\{h = E\}}
\]
is also a $K$-invariant contact form.  By the implicit function
theorem, for $E>0$ sufficiently small, there is a function $f: \{q =
E\} \to (0, \infty)$, which is $C^2$ close to 1, so that
\[
\phi :\{ q = E\} \to \{h = E\} \quad \phi (x) = f(x) x
\]
is a $K$-equivariant diffeomorphism.  Since $\phi ^* \beta _E = f^2
\alpha _E$, the manifolds $\{q = E\}$ and $\{h = E\}$ are
$K$-equivariantly contactomorphic. Moreover, under the identification $\phi$
the two contact forms $\alpha_E$ and $\beta_E$ are $C^2$-close (again,
provided $E$ is small).  Note that the two associated contact moment maps are 
$\Phi |_{\{q = E\}}$ and $\Phi |_{\{h= E\}}$ respectively.\\

Up to re-parameterization the integral curves of the Hamiltonian
vector field of $h$ in $\{h =E\}$ are the integral curves of the Reeb
vector field $X_E$ of $\beta_E$.  Similarly the integral curves of the
Hamiltonian vector field of $q$ on $\{q = E\}$ are the integral curves
of the Reeb vector field $Y_E$ of $\alpha_E$.  In particular the
relatively periodic orbits of $h$ on $\{h = E\}$ are relatively
periodic orbits of $X_E$.  Since the hypersurface $\{h = E\}$ is
compact and since the orbit type decomposition of the contact quotient
$\{h = E\}/\!/ G$ is a stratification, the minimal strata of the
quotient are compact.  Let $Q = (\{h = E\}/\!/ G)_{(L)}$ be one such
stratum.  Then the relatively periodic orbits of $h$ in $\{h = E\}
\cap \Phi \inv (0) \cap V_{(L)}$ descend to periodic orbits of the
Reeb vector field $X$ of the contact form $\beta _{(L)}$ on $Q = (\{h
= E\} \cap \Phi \inv (0) \cap V_{(L)})/G$.  Therefore to prove
Theorem~\ref{theorem2} it is enough to establish the existence of
periodic orbits of $X$.  For this, according to Theorem~\ref{theorem3}, it
suffices to establish the existence of a contact form $\alpha $ on $Q$
whose Reeb vector field $Y$ generates a torus action and such that 
$\beta _{(L)}$ is $C^2$ close to $\alpha$ when $E>0$ is small enough.

The form $\alpha$, of course, is the form induced by $\alpha_E$.  Let
us prove that it does have the desired properties. Since $\phi : \{q =
E\}\to \{h = E\}$ is an equivariant contactomorphism it induces an
identification of the contact manifold $Q$ with $ (\{q = E\} \cap \Phi
\inv (0) \cap V_{(L)})/G$.  Moreover, since $\alpha_E$ and $\beta_E$
are $C^2$-close, the induced forms $\beta_{(L)}$ and $\alpha = \alpha
_{(L)}$ are $C^2$-close as well.  Since $q$ is definite, its
Hamiltonian flow generates a linear symplectic action of a torus $\TT$
on $V$.  The restriction of this action to $\{ q = E\}$ is also
generated by the Reeb vector field $Y_E$ of $\alpha _E$.  Since $q$ is
$K$-invariant, the action of $\TT$ commutes with the action of $K$ and
preserves the moment map $\Phi$.  Hence it descends to an action of
$\TT$ on $Q$.  Moreover, since the Reeb vector field of $\alpha_E$
descends to the Reeb vector field of $\alpha _{(L)}$ on $Q$, the
induced action of $\TT$ on $Q$ is generated by the Reeb vector field of
$\alpha _{(L)}$.  We conclude that Theorem~\ref{theorem3} implies
Theorem~\ref{theorem2}.

\section{ Perturbations of Reeb flows: proof of Theorem~\ref{theorem3}}

In the proof of Theorem~\ref{theorem3} we will need the following 
elementary result.
\begin{lemma} \label{last-lemma}
Let $\phi_t$ be a dense one-parameter subgroup in a torus $\TT$ and let
$H$ be a subgroup of $\TT$ topologically generated by an element
$\phi^\tau$, $\tau>0$.  Then either $H$ has codimension one in $\TT$ or 
$H=\TT$.
\end{lemma}

\begin{proof}[Proof of Lemma~\ref{last-lemma}]
It suffices to show that the map
$$
[0,\tau]\times H\to \TT, \quad F(t,h)=\phi^t\cdot h
$$
is onto $\TT$.
Pick $g\in \TT$. Assume first that $g$ is in the one-parameter subgroup, 
i.e., $g=\phi^t$ for some $t$. Then we have
$t=k\tau+t'$ with $0\leq t' < \tau$ and, clearly,
$$
g = \phi^{t'}\cdot [(\phi^{\tau})^k] = F\left(t', (\phi^\tau)^k\right).
$$
Hence $g$ is in the image of $F$.

Let now $g$ be in $\TT$, but not in the one-parameter subgroup $\phi^t$.
Then there exists a sequence
$t_r\to\pm\infty$ such that $\phi^{t_r}\to g$. (This sequence must go
to positive or negative infinity, for otherwise $g$ is in the one-parameter
subgroup.) Assume that $t_r\to\infty$; the case of negative infinity can be
dealt with in a similar fashion. As above, we write
$$
t_r= k_r \tau + t'_r,
$$
where $k_r\to \infty$ as $r\to \infty$ and $0 < t'_r <\tau$.

The elements $(\phi^\tau)^{k_r}$ are in $H$ and, since $H$ is compact,
we may assume that  $(\phi^\tau)^{k_r}\to h\in H$ by passing if
necessary to a subsequence. Furthermore, by passing if necessary to a
subsequence again, we may assume that $t'_r\to t'\in [0,\tau]$.

We claim now that $g=F(t',h)$. To see this note that as above
$$
\phi^{t_r}=\phi^{k_r \tau + t'_r}=\phi^{t'_r}\cdot [(\phi^\tau)^{k_r}].
$$
As $r$ goes to infinity, $\phi^{t'_r}\to\phi^{t'}$ and the
second term goes to $h$. Hence,
$$
g=\phi^{t'}\cdot h=F(t',h).
$$
This completes the proof of the lemma.
\end{proof}

\begin{proof}[Proof of Theorem~\ref{theorem3}]
First, let us set notation. We denote by $X$ the Reeb
vector field of $\alpha$ and by $\phi^t$ its Reeb flow.  By the
hypotheses of the theorem, the flow $\phi^t$ generates an action of a
torus $\TT$ on $Q$. We will view $\phi^t$ as a dense one-parameter
subgroup of $\TT$.  The points on periodic orbits of $X$ will be
referred to as periodic points.
We break up the proof of the theorem into four steps. Steps 1--3 concern 
exclusively properties of the Reeb flow of $\alpha$. The perturbed form $\beta$
enters the proof only at the last step.

1. We claim that the periodic points of $X$ are exactly the points $x\in Q$
whose stabilizers $\TT_x$ have codimension one in $\TT$.

Indeed, let $x\in Q$ be a periodic point, i.e., $\phi^T(x)=x$ for some
$T>0$.  Since $\phi^t$ is dense in $\TT$, the Reeb orbit through $x$ is
dense in the $\TT$-orbit through $x$. Since $x$ is a periodic point, the
Reeb orbit is closed and thus equal to the $\TT$-orbit. Hence, $\TT/\TT_x$
is a circle and thus $\TT_x$ has codimension one. Conversely, if $\TT_x$
has codimension one, the Reeb orbit through $x$ must be dense in the
$\TT$-orbit and hence equal to the $\TT$-orbit because the latter is a
circle.

2. Let now $N$ be a minimal stratum of the $\TT$-action, which is
comprised entirely of periodic points. We claim that such a stratum
exists, is a smooth submanifold, and all points of $N$ have the same
period, i.e., the $\TT$-action on $N$ factors through a free circle
action.

Since the $\TT$-action has no fixed points, periodic points lie in
minimal strata of the action. Furthermore, the Reeb flow of
$\alpha$ has at least one periodic orbit (in fact, at least two unless
$Q$ is a circle); this
follows, for example, from a theorem of Banyaga and Rukimbira,
\cite{BR}. Now it suffices to take as $N$ a minimal stratum containing a 
periodic point.  The fact that $N$ is smooth is a general
result about compact group actions. Finally, all points in $N$ have
the same stabilizer $\TT_x$ and the action of the circle $\TT/\TT_x$ on $N$ is
free because $N$ is minimal. The period $T$ of $x\in N$ is the first 
$T>0$ such that $\phi^T \in \TT_x$.

3. We claim that $N$ is a non-degenerate invariant submanifold for the
Reeb flow of $\alpha$. 

Let $x\in N$. We need to show that the linearization $d\phi^T$ on the
normal space $\nu_x$ to $N$ at $x$ does not have unit as an
eigenvalue. By definition, this linearization is just the linearized
action of $\phi^T \in \TT_x$ on $\nu_x$. As is well known, the
isotropy representation of $\TT_x$ on $\nu_x$ contains no trivial
representations in its decomposition into the sum of irreducible
representations. Hence, it suffices to show that the subgroup
generated by $\phi^T$ is dense in $\TT_x$.

Let $m$ be the first positive integer such that $\left(\phi^T\right)^m$ 
is in $\TT_x^0$, the
connected component of identity in $\TT_x$. (Such an integer $m$ exists because
$\TT_x/\TT_x^0$ is a finite subgroup of
the circle $\TT/\TT_x^0$.) Since $\TT_x/\TT_x^0$ is finite cyclic, 
it suffices to show 
that the subgroup
$H$ topologically generated by $\phi^{mT}$ is equal to $\TT_x^0$. This 
follows immediately from Lemma~\ref{last-lemma}. Indeed, by the lemma,
the group $H$ is either equal to $\TT$ or
has codimension one in $\TT$. Since $H\subset \TT_x^0$ and $\TT_x^0$ has
codimension one, the group $H$ must have codimension one. Thus $H$ is a
closed subgroup of $\TT_x^0$ of the same dimension as $\TT_x^0$ and hence
$H=\TT_x^0$.

4. Now we invoke the following theorem due to Kerman \cite{Ke}
(p. 967). Let $Crit(P)$ be the minimal possible number of critical
points of a smooth function on a compact manifold $P$.

\begin{theorem}[Kerman, \cite{Ke}] 
Let $Q$ be a compact odd-dimensional manifold, $X$ a
non-vanishing vector field on $Q$, and $N$ a non-degenerate periodic 
submanifold of $X$. Let $\Omega$ be a closed maximally non-degenerate
two-form on $Q$ whose kernel is $C^2$-close to $X$ and such that the
class $[\Omega|_N]$ is in the image of the pull-back from $H^2(N/S^1)$ to 
$H^2(N)$. Then $\Omega$ has at least $Crit(N/S^1)$ closed characteristics
near $N$.
\end{theorem}

Applying this theorem to $Q$, $N$ and $X$ as above, and $\Omega=d\beta$
we obtain the required result.                            
\end{proof}

\begin{remark}
In fact, our proof of Theorem~\ref{theorem3} establishes the existence of
two distinct periodic orbits when $Q$ is not a circle. As a consequence, in 
the setting of Theorems~\ref{main} and \ref{theorem2} there exist at least 
two distinct relative periodic orbits unless $Q$ is a circle.
\end{remark}


\end{document}